\pdfoutput=1
\RequirePackage{ifpdf}
\ifpdf 
\documentclass[pdftex]{sigma}
\else
\documentclass{sigma}
\fi

\numberwithin{equation}{section}

\begin{document}

\allowdisplaybreaks

\renewcommand{\PaperNumber}{020}

\FirstPageHeading

\ShortArticleName{On a~Seminal Paper by Karlin and McGregor}

\ArticleName{On a~Seminal Paper by Karlin and McGregor}

\Author{Mirta M.~CASTRO~$^\dag$ and F.~Alberto GR\"UNBAUM~$^\ddag$}

\AuthorNameForHeading{M.M.~Castro and F.A.~Gr\"unbaum }

\Address{$^\dag$~Departamento Matem\'atica Aplicada II, Universidad de Sevilla,
\\
\hphantom{$^\dag$}~c$\backslash$Virgen de \'Africa 7, 41011, Sevilla, Spain}
\EmailD{\href{mailto:mirta@us.es}{mirta@us.es}}

\Address{$^\ddag$~Department of Mathematics, University of California, Berkeley, Berkeley, CA 94720 USA}
\EmailD{\href{mailto:grunbaum@math.berkeley.edu}{grunbaum@math.berkeley.edu}}

\ArticleDates{Received August 04, 2012, in f\/inal form February 25, 2013; Published online March 02, 2013}

\Abstract{The use of spectral methods to study
birth-and-death processes was pioneered by S.~Karlin and
J.~McGregor.
Their expression for the transition probabilities was made explicit
by them in a~few cases.
Here we complete their analysis and indicate a~few applications
of their very powerful method.}


\Classification{15A24; 39A99; 47B36; 60J80}

\section{Introduction}
\subsection{Preliminaries}

In a~paper overf\/lowing with new as well as classical ideas, S.~Karlin and
J.~McGregor~\cite{KMcG} considered
the birth-and-death process whose one step transition
probability matrix ${\mathbb P}$ is given by
\begin{gather*}
{\mathbb P}=\begin{pmatrix}
r_0&p_0&0&0
\\
q_1&r_1&p_1&0
\\
0&q_2&r_2&p_2
\\
&&\ddots&\ddots&\ddots
\end{pmatrix},
\end{gather*}
where all indicated entries are nonnegative, $p_0+r_0\leq1$ and $p_i+r_i+q_i\leq1$, $i\geq1$.
These inequalities mean that there is a~positive probability of jump to an (ignored) cof\/f\/in
state, or even distinct cof\/f\/in states for each $i$ such that the corresponding sum is not~1.

The paper~\cite{KMcG} acknowledges and builds on two dif\/ferent earlier approaches: one
consists in exploiting the spectral theory of an operator in an appropriate Hilbert space,
see~\cite{F1,K,McK}
and the other one relies on a~hands-on study of sample path properties, see~\cite{F2, H,HR}.

They introduce the polynomials $Q_j(x)$ by the conditions
$Q_{-1}(0)=0$, $Q_0(x)=1$ and using the notation
\begin{gather*}
Q(x)=\begin{pmatrix}
Q_0(x)
\\
Q_1(x)
\\
\vdots
\end{pmatrix}
\end{gather*}
they enforce the recursion relation
\begin{gather*}
{\mathbb P}Q(x)=x Q(x)
\end{gather*}
to argue the existence of a~unique measure $\psi(dx)$ supported in $[-1,1]$ such that
\begin{gather*}
\pi_j\int_{-1}^1Q_i(x)Q_j(x)\psi(dx)=\delta_{ij}.
\end{gather*}
Then they obtain what is nowadays called the Karlin--McGregor representation formula
\begin{gather*}
\big({\mathbb P}^n\big)_{ij}=\pi_j\int_{-1}^1x^nQ_i(x)Q_j(x)\psi(dx),
\end{gather*}
where the constants $\pi_j$ are given by
\begin{gather}
\label{normas}
\pi_0=1,
\qquad
\pi_j=\frac{p_0p_1\cdots p_{j-1}}{q_1q_2\cdots
q_j},
\qquad
j\ge1.
\end{gather}
Many general results can be obtained from the representation formula
given above; as they remark, this is nothing but the spectral theorem.

After establishing many such results the authors turn their attention to the study of the spectral
properties in the case when one has essentially constant transition probabilities, more
specif\/ically
\begin{gather}\label{matriztridi}
{\mathbb P}=\begin{pmatrix}
r_0&p_0&0&0
\\
q&0&p&0
\\
0&q&0&p
\\
&&\ddots&\ddots&\ddots
\end{pmatrix}.
\end{gather}
Here the conditions on the parameters $r_0$, $p_0$, $q$ and $p$ are
\begin{gather*}
r_0\ge0,\qquad p_0>0,
\qquad
r_0+p_0\le1,
\qquad
p,q>0,
\qquad
p+
q=1.
\end{gather*}
One usually puts $q_0=1-r_0-p_0$ and when $q_0>0$ there is a~loss of total mass of the
initial probability distribution.

If one denotes the spectral measure by $\psi(dx)$ and by $m(z)$ its
Stieltjes transform
\begin{gather*}
m(z)=\int_{-1}^1\frac{\psi(dx)}{x-z},
\end{gather*}
there is a~classical method to obtain $m(z)$, which we recall now.
If $n(z)$ denotes
the Stieltjes transform of the measure going along with the ``chopped matrix'' obtained
by deleting the f\/irst row and column of ${\mathbb P}$ then one proves that
\begin{gather*}
m(z)=-1/(z-r_0+q p_0n(z)).
\end{gather*}
The computation of $n(z)$ itself can be made using the same trick; moreover since (in the case of
an essentially constant
transition matrix) the ``twice
chopped matrix'' is the same as the ``chopped matrix'' we get that $n(z)$ satisf\/ies $
n(z)=-1/(z+q p n(z))$
or $n(z)^2p q+z n(z)+1=0$, which gives, by solving this quadratic equation
and observing that $n(z)$ should vanish at inf\/inity,
$n(z)=\big(-z+\sqrt{z^2-4p q}\big)/(2p q)$.

Before proceeding, it is useful to note that the function $n(z)$ is well def\/ined as
a single valued function in the complex plane from which we have removed the closed
interval $[-2\sqrt{pq},2\sqrt{pq}]$.
As one approaches this cut from above, the
function $n(z)$ has a~nontrivial imaginary part coming directly from the square root.
This square root has positive values on the real axis to the right of the cut and negative
values on the real axis to the left of the cut.
This basic observation will play an important role
at several points below.

Returning to the expression above one gets for $m(z)$
\begin{gather*}
m(z)=-\frac{1}{z-r_0+ \dfrac{p_0}{2p}\bigl({-}z+\sqrt{z^2-4p q} \bigr)}.
\end{gather*}
This expression may have zeros in the denominator giving rise to poles in $m(z)$.
Observe that none of these zeros can take place for values of $z=x$ in the cut itself.
This is, for instance, a~consequence of the fact that along the cut the square root term gives a
nontrivial imaginary contribution that cannot be cancelled by the other terms in the
denominator which are real valued.

A good way to study this function $m(z)$ is to get rid of the square root in the denominator.
In this fashion one gets the expression for~$m(z)$
already given in
\cite[Section 2ii]{KMcG}, namely
\begin{gather}
\label{m(z)}
 m(z)=\frac{r_0-\left(1-\dfrac{p_0}{2p}\right)z+\dfrac{p_0}{2p}\sqrt{z^2-4pq}}{\left(1-
 \dfrac{p_0}{p}\right)z^2-2r_0\left(1- \dfrac{p_0}{2p}
\right)z+r_0^2+ \dfrac{p_0^2q}{p}}.
\end{gather}
There is a~classical way of recovering $\psi(dx)$ from $m(z)$ which
consists of considering~$m(z)$ as meromorphic function def\/ined on an
appropriate Riemann surface and inverting this Stieltjes transform.
The measure consists of a~nice density
plus possibly some delta masses with certain weights.
The function~$m(z)$
is well def\/ined in the same region that makes~$n(z)$ single-valued.
The continuous part of the measure results from
looking at the values of the imaginary part of~$m(z)$ on one side of the cut.
The discrete masses in~$\psi(dx)$ come from the residues at the simple poles of~$m(z)$ but a~more
detailed analysis is given below.
Clearly this inversion process requires a~bit of careful computation.

There are several references in the literature where the explicit computation of distributions
associated to some special classes of Jacobi matrices is handled.
One can mention~\cite{Ge2} (see also~\cite{Apt}) for the case of periodic
coef\/f\/icients,~\cite{Al1, Al} for the case of parameters with constant limits,~\cite{DKS} for
asymptotically periodic parameters and~\cite{AW} in a~dif\/ferent general context.

For many other aspects of birth-and-death processes the reader can see~\cite{Ch,F,G,ILMV,K1,LR}.
Many of these reference deal with continuous time.

\subsection{The story behind this paper}

In~\cite{KMcG} one f\/inds two examples done in some detail.
Specif\/ically they deal with the case when $r_0=0$, and the case when $p_0=p$, $r_0=q$, which insures
that the total mass of the initial probability
distribution is preserved.
For $r_0=0$ and $p=p_0$, which gives rise to a~variant of the Wigner semicircle law, the authors
also consider a~two-periodic transition probability matrix.

As preparation for a~larger project we had planned to complete the
analysis given in~\cite{KMcG} by considering the case given in~\eqref{matriztridi}.
Fortunately, during this process, Michael Anshelevich pointed us
in the direction of some references in a~dif\/ferent f\/ield, and after appropriate conversions it
became clear that a~full analyis of the problem posed in~\cite{KMcG} was possible
by using results in the literature.

The most recent relevant paper is~\cite{SY} in the context of ``free
probability'' and the reader may also want to consult~\cite{A}.
Our
basic reference is an earlier paper~\cite{CT} where Stieltjes
inversion along the lines indicated above is avoided altogether.
The
method in~\cite{CT} consists of using an old method of Christof\/fel
which can be found in the very classical reference~\cite{Sz}, and then
getting a~nice corollary which amounts to using the method of
Christof\/fel in reverse.

In~\cite{SY} the authors observe that there is a~small error in
\cite{CT} amounting to having normalized the continuous and the
discrete parts of $\psi(dx)$ in dif\/ferent fashions.
The method used
in~\cite{SY} is the method of Stieltjes indicated above.
In fairness
one should say that in~\cite{SY} there are several details left out
but that the f\/inal result, essentially already given in~\cite{CT}, is
correct.
See also~\cite{A1}.

One should remark that neither of the references~\cite{CT, SY}
mention the results in~\cite{KMcG}.
A look at more recent papers indicates
that the results in~\cite{KMcG} are not widely known and that several of
these cases are being rederived separately.
For instance,~\cite{S} does not
refer
to any of the papers~\cite{CT,KMcG, SY} although the author deals with
several special cases of the problem considered in these references, such as the
Wigner, Marchenko--Pastur and Kesten--McKay measures.

Going back to the problem considered in~\cite{KMcG} and by appropriately adapting the results
in~\cite{CT, SY} one can see
that the general expression for the measure $\psi(dx)$ is given as
follows.
The continuous part of the measure is given by
\begin{gather*}
\frac{\dfrac{p_0}{2\pi} \sqrt{4pq-x^2}}{(p-p_0)
x^2-r_0(2p-p_0)x+pr_0^2+p_0^2q} dx,
\qquad
|x|\le 2\sqrt{pq}.
\end{gather*}
To this continuous part one needs to add possibly one
or two mass points, and at this point the analysis in~\cite{CT,SY} becomes relevant.

In this paper we do not reproduce the arguments given in~\cite{CT,SY} to obtain the orthogonality
measure.
Instead we display the results in a~``reader friendly" way, we go through the steps showing
why certain residues turn out to be positive, and more importantly we exploit the spectral method
discussed here to derive
some information about ``quasi-stationary distributions".
This goes beyond the analysis in~\cite{CT,KMcG,SY}
and combines this analysis with important work by E.~van Doorn and P.~Schrijner, see for instance~\cite{vDSc,vDSc1}.
Some of these results may be hard to obtain unless one uses spectral methods as in this paper
and gives a~rather
nice way to complete the analysis started in~\cite{KMcG} a~long time ago.

Since our expressions dif\/fer slightly from those in~\cite{CT} and \cite{SY} we have checked
carefully the internal consistency of our choices, including all normalizations.
This required determining the eigenfunctions of the matrix~\eqref{matriztridi}.
Eventually we found out that~\cite[p.
204]{Ch} referring to earlier work of Geronimus~\cite{Ge1} had already essentially written down
these eigenfunctions as
\begin{gather}
\label{polynoms}
Q_j(x)=\left(\frac{q}{p}\right)^{\frac{j}{2}}
\left[\frac{2(p_0-p)}{p_0} T_j (x^* )+\frac{2p-p_0}{p_0} U_j(x^*)-\frac{r_0}{p_0} \sqrt{\frac{p}{q}} U_{j-1}(x^*)\right],
\end{gather}
where $T_j$ and $U_j$ are the Chebyshev polynomials of the f\/irst and
second kind, and $x^*=x/(2\sqrt{pq})$.

When the measure is made up of the continuous part given above plus the
po\-ssible contributions coming from the residues at the simple poles of $m(z)$, the squared norms
of these polynomials are as in~\eqref{normas}:
\begin{gather*}
\pi_0=1,
\qquad
\pi_j=\frac{p_0p^{j-1}}{q^j},
\qquad
j\geq1.
\end{gather*}

\section{Beyond Karlin and MacGregor}

\subsection{A single pole}
\label{una raiz}

{\sloppy The denominator of $m(z)$ has a~single root exactly when $p_0=p$, $r_0\neq0$, in which
case the root is
\begin{gather}
\label{single root}
x_1=r_0+p q/r_0.
\end{gather}
Here is a~simple argument showing that the root has to be to the right of the
interval $[-2\sqrt{pq},2\sqrt{pq}]$.
Indeed, we must have
that $r_0+p q/r_0$ should be at least $2\sqrt{pq}$, and the
dif\/ference of these two quantities can be rewritten as $\big(\sqrt{r_0}-\sqrt{pq/r_0} \big)^2.$

}

In the boundary case $r_0=\sqrt{pq}$ the simple root is located at $2\sqrt{pq}$ and it will turn
out that we have no mass.

For $m(z)$ we get the expression
\begin{gather*}
m(z)=\frac{1}{2r_0} \frac{z-2r_0-\sqrt{z^2-4pq}}{z-(r_0+pq/r_0)}.
\end{gather*}
We are interested in the residue of $m(z)$ at~\eqref{single root} and this calls for a~careful
evaluation of the
numerator at $z=x_1.$

\begin{figure}[t]
\centering
\begin{minipage}{75mm}
\includegraphics[width=70mm]{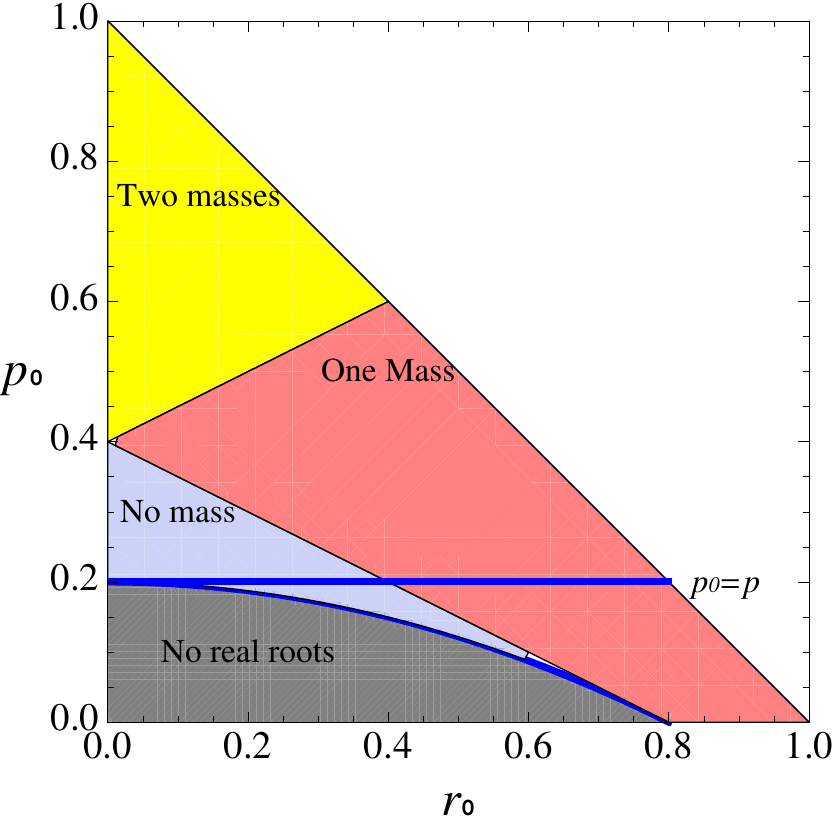}
\caption{Case $p=0.2$.}
\label{Caso p0.2}
\end{minipage}\qquad
\begin{minipage}{75mm}
\includegraphics[width=70mm]{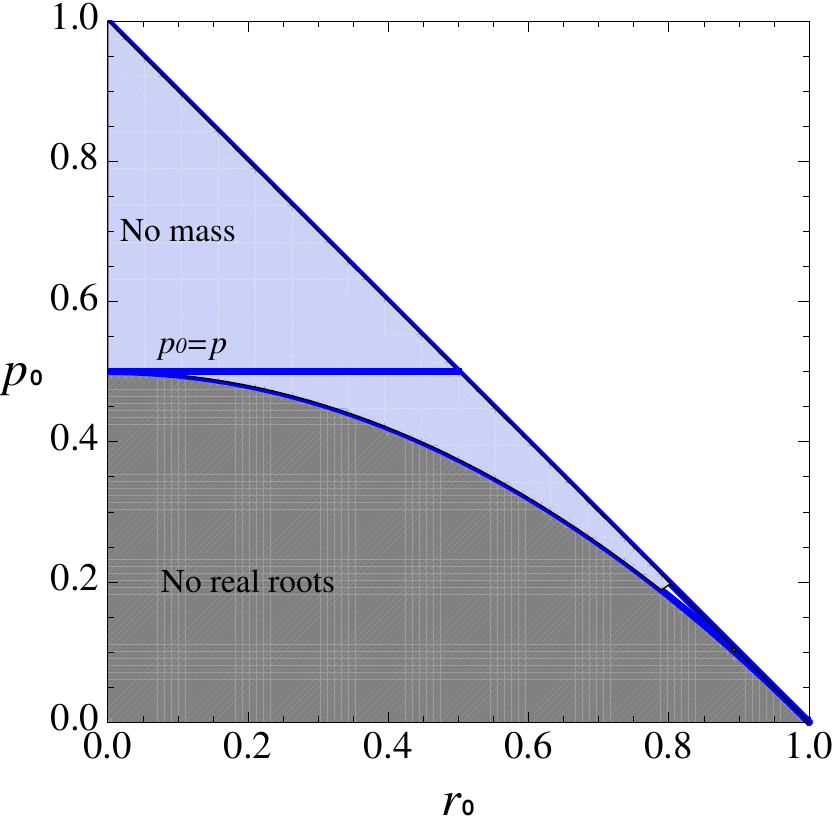}
\caption{Case $p=0.5$.}
\label{Caso p=0.5}
\end{minipage}
\end{figure}

\begin{figure}[t]

\centering
\includegraphics[width=70mm]{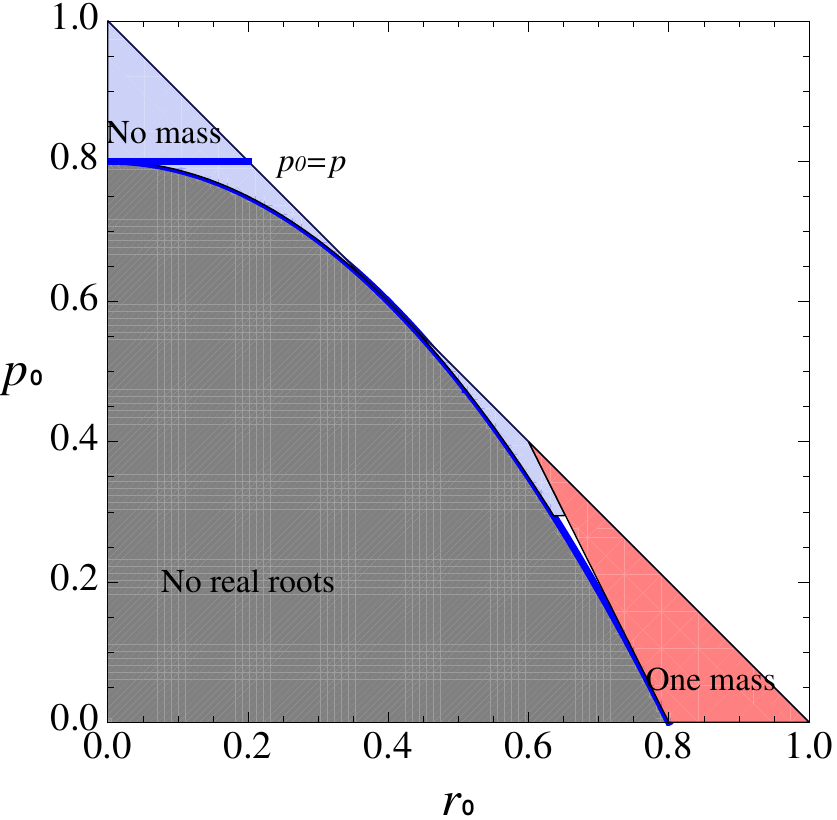}
\caption{Case $p=0.85$.}
\label{Caso p=0.85}
\end{figure}

Under the square root we have $
(r_0+pq/r_0)^2-4p q=(r_0-pq/r_0)^2.
$ Since we are using the determination of the square root function which is positive on the real
axis to the right of $2\sqrt{pq}$ and we have shown that $x_1\geq2\sqrt{pq}$ we have for the
residue the value
\begin{gather*}
\frac{1}{2r_0} \left(r_0+\frac{pq}{r_0}-2r_0-\left|r_0-\frac{pq}{r_0}\right|\right).
\end{gather*}
This can be rewritten as $
(1/2)(-(1-pq/r_0^2)-|1-pq/r_0^2|),
$ which gives $-(1-pq/r_0^2)_+$, where the subindex $+$ denotes, as usual, the positive part.
Since the denominator in the def\/inition of~$m(z)$ is $x-z$ instead of $z-x$ the weight attached
to this root is given by the negative of the residue, namely
\begin{gather}
\label{pos part p=p0}
\left(1-\frac{pq}{r_0^2}\right)_+.
\end{gather}
Here and later it is important to note that the presence of zeroes in
the denominator of $m(z)$ does not guarantee the presence of a~mass
at such a~zero, since the residue may vanish.
But if this residue does not vanish then it is automatically positive as we have seen in the simple
case analyzed above.
In order to have a~mass when $p_0=p$ one should have
\begin{gather}
\label{una masa p=p0}
p<1/2
\qquad
\textrm{and}
\qquad
r_0\in\big(\sqrt{pq},1-p\big],
\end{gather}
see Figs.~\ref{Caso p0.2},~\ref{Caso p=0.5} and~\ref{Caso p=0.85}.

\subsection{Two roots}\label{dos raices}

In the case when the denominator has two dif\/ferent zeros at locations
$x_i$, $i=1,2$, an analysis of the residues along the lines of the simpler case above gives for the
values of the weights
the expressions
\begin{gather}
\label{pos part}
\frac{1}{\sqrt{r_0^2-4q(p-p_0)}}\left(\frac{qp_0^2}{p|x_i-r_0|}-|x_i-r_0|\right)_+.
\end{gather}
In particular, we see that to have a~mass at $x_i$ we need
\begin{gather}
\label{cond masa 1}
r_0^2\ge4q(p-p_0)
\end{gather}
as well as
\begin{gather}
\label{condicion masa 2}
(x_i-r_0)^2<\frac{q p_0^2}{p}.
\end{gather}

{\samepage
We will describe, in terms of the values of the parameters $r_0$, $p_0$ and $p$, the dif\/ferent
possibilities according to the number of positive masses we may have.
We will show below that we have just one mass exactly when
\begin{gather}
p<\frac12
\qquad
\textrm{and}
\qquad
2p- \sqrt{\frac{p}{q}}r_0<p_0\leq2p+ \sqrt{\frac{p}{q}} r_0, \qquad \textrm{or}\label{cond one mass 1}
\\
p>\frac12
\qquad
\textrm{and}
\qquad
p_0>2p- \sqrt{\frac{p}{q}} r_0,
\qquad
r_0>\sqrt{pq},
\label{con one mass 2}
\end{gather}
see Figs.~\ref{Caso p0.2} and \ref{Caso p=0.85}.
Here, the mass point is always located on the right of the cut.
}

In the same way, we have two masses exactly when
\begin{gather}
\label{dos masas}
p<\frac 12
\qquad
\textrm{and}
\qquad
p_0>2p+ \sqrt{\frac{p}{q}} r_0,
\end{gather}
see in particular Fig.~\ref{Caso p0.2}.
In this case the two mass points are located at both sides of the cut in
$[-1,1]\setminus[-2\sqrt{pq},2\sqrt{pq}]$.

Indeed, we write the explicit expression of each $x_i$, $i=1,2$:
\begin{gather*}
x_{1,2}=\frac{r_0(2p-p_0)\pm p_0\sqrt{r_0^2-4q(p-p_0)}}{2(p-p_0)},
\qquad
p_0\neq p.
\end{gather*}
In order to have a~mass at $x_i$, $i=1,2$, considering~\eqref{condicion masa 2}, one needs to have
\begin{gather}
\label{desig masas}
\left|r_0\pm\sqrt{r_0^2-4q(p-p_0)}\,\right|<2\frac{\sqrt{q}}{\sqrt{p}}|p-p_0|.
\end{gather}
To analyze carefully inequality~\eqref{desig masas} one has to consider separately the cases
$p_0<p$ and $p_0>p$.

Suppose that $p_0<p$.
We will show that in this case it is possible to have a~mass only at $x_2$ when
\begin{gather}
\label{desmenos}
p_0>2p- \sqrt{\frac{p}{q}} r_0
\end{gather}
holds true.
Indeed, if we suppose that there is a~positive mass at $x_1$, from~\eqref{desig masas} we can write
\begin{gather*}
r_0+\sqrt{r_0^2-4q(p-p_0)}<2\frac{\sqrt{q}}{\sqrt{p}}(p-p_0),
\end{gather*}
thus one should have
\begin{gather*}
2\frac{\sqrt{q}}{\sqrt{p}}(p-p_0)>r_0,
\end{gather*}
that contradicts~\eqref{cond masa 1}.
Indeed, considering the previous inequality and then~\eqref{cond masa 1} one would have{\samepage
\begin{gather*}
r_0^2<\frac{4q}{p}(p-p_0)^2\leq\left(1-\frac{p_0}{p}\right)r_0^2,
\end{gather*}
which is a~contradiction.}

Analogously, to have a~positive mass at $x_2$ we write for~\eqref{desig masas}
\begin{gather*}
  r_0-\sqrt{r_0^2-4q(p-p_0)}<2\frac{\sqrt{q}}{\sqrt{p}}(p-p_0),
\end{gather*}
or equivalently
\begin{gather*}
r_0-2\frac{\sqrt{q}}{\sqrt{p}}(p-p_0)<\sqrt{r_0^2-4q(p-p_0)}.
\end{gather*}
If~\eqref{desmenos} holds true then the expressions in both sides of the previous inequality are
nonnegative so we can consider the squared inequality, which happens to be equivalent
to~\eqref{desmenos} as well.
We remark that when $p_0<p$~\eqref{desmenos} implies $r_0>\sqrt{pq}$.

We conclude the analysis of the case $p_0<p$ showing that if $x_2$ has a~positive mass then
$x_2>2\sqrt{pq}$.
Indeed, the inequality $x_2-2\sqrt{pq}>0$ is equivalent to
\begin{gather*}
r_0(2p-p_0)-4\sqrt{pq}(p-p_0)>p_0\sqrt{r_0^2-4q(p-p_0)}.
\end{gather*}
Using~\eqref{cond masa 1} one verif\/ies that the left part of the previous inequality is
nonnegative thus taking the squares and simplifying conveniently one obtains the equivalent
expression
\begin{gather*}
\big(r_0\sqrt{p}-\sqrt{q}(2p-p_0)\big)^2>0.
\end{gather*}
We point out that in the boundary case $p_0=2p- \sqrt{p/q} r_0$ the weight for $x_2=2\sqrt{pq}$
in~\eqref{pos part} is zero.

For the case $p_0>p$, using similar arguments as before, one can see that to have a~mass at~$x_1$
one needs to have
\begin{gather}
\label{desmas}
p_0>2p+ \sqrt{\frac{p}{q}} r_0.
\end{gather}
Since
\begin{gather}
\label{suma1}
p_0+r_0\leq1,
\end{gather}
condition~\eqref{desmas} implies $p<1/2$.
Analogously, if~\eqref{desmenos} holds true, there is a~positive mass at~$x_2$.
Consequently, in this case one has two masses when~\eqref{dos masas} holds true, and one verif\/ies
that $x_1<-2\sqrt{pq}$ and $x_2>2\sqrt{pq}$.

Notice that, in both cases described above, in order to have a~mass one should have $p\neq1/2$,
otherwise we will get a~contradiction with condition~\eqref{suma1} since the parameters must
verify~\eqref{desmenos}.
See Figs.~\ref{Caso p0.2},~\ref{Caso p=0.5} and \ref{Caso p=0.85}.

One can consider the case when our process is recurrent (see~Sections~2 in \cite{KMcG,vDSc,vDSc1})
for characterizations of some probabilistic properties for this kind of model), which holds when
\begin{gather}
\label{cond recurrencia}
r_0+p_0=1
\qquad
\textrm{and}
\qquad
p\leq \tfrac{1}{2}.
\end{gather}
In this case $z=1$ is a~root of the denominator in $m(z)$, and has a~mass here only when the
process is positive recurrent as noticed in~\cite[Section~2]{KMcG}, which corresponds to $p<1/2$
in~\eqref{cond recurrencia}.
As a~very special subcase, one considers $r_0=0$, $p_0=1$, corresponding to a~ref\/lecting boundary
condition at the left most
state.
Then one has roots of $m(z)$ both at $z=1$ and $z=-1$, having weights only when $p<1/2$, as pointed
out initially in~\cite[Section~2]{KMcG} and \cite[pp.~108--109]{K1}.

\subsection{Some useful f\/igures}

Here we have included some useful f\/igures illustrating the dif\/ferent situations that may arise for
the roots of the denominator of~$m(z)$ in~\eqref{m(z)}, discussed above.
The nature of the f\/igures depends on the value of a~f\/ixed value of the parameter $p$ being
less, equal or larger than~$1/2$.
We consider free nonnegative parameters $p_0$ and $r_0$ satisfying~\eqref{suma1}.
As noticed before, in the case of Fig.~\ref{Caso p0.2} for $p<1/2$ one can see that there may
be one or two mass points.
These situations are delimited by the lines $p_0=2p- \sqrt{p/q} r_0$ and $p_0=2p+ \sqrt{p/q} r_0$.
We point out that in order to have two masses it is necessary that $p_0>p$.

For $p=1/2$ Fig.~\ref{Caso p=0.5} shows that we do not have any mass points.
For $p>1/2$, as Fig.~\ref{Caso p=0.85} shows, there may be one mass point above the line
$p_0=2p- \sqrt{p/q} r_0$ and for $p_0<p$.

\section{Miscellaneous}

\subsection{Quasi-stationary distribution}

A new notion enters at this point.
A quasi-stationary distribution for a~Markov chain with an absorbing state (the case of our
chain when~$q_0$ does not vanish) is def\/ined as an initial distribution such that the
(conditional) probability of the chain
being at state~$j$ at time~$n$, conditioned on the fact that no absorption has occurred by then, is
independent of the time~$n$ for all states~$j$.

The issue of the existence and/or uniqueness of such distributions can be illuminated by the
analysis given above.

Assume that $q_0>0$ so that we have a~loss of total mass for any initial probability distribution.
Assume further that
$p\leq1/2$ so that eventual transition out of the state space is guaranteed.
We denote by $p_{-1}(n)$ the probability of having left the state space by time~$n$.

One def\/ines a~(normalized) probability distribution ${\alpha_i}$ on the nonnegative integers as
a~quasi-stationary distribution for our walk if, by denoting with $p_j(n)$ the probability that the
walk be at state j at time n for the random walk with initial distribution $p_j(0)=\alpha_j$, one
gets the relation
\begin{gather*}
\frac{p_j(n)}{1-p_{-1}(n)} =\alpha_j,
\qquad
\textrm{for}
\quad
j=0,1,2,\ldots,
\quad
n=0,1,2,\ldots.
\end{gather*}

Denote now by $\eta$ the supremum of all points in the support of the orthogonality measure.
This can be a~point mass to the right of $2\sqrt{pq}$ or the value
$2\sqrt{pq}$ itself.
In both cases one can have $\eta=1$.
This value plays a~very important role in terms of the existence and uniqueness of quasi-stationary
distributions as follows from a~more general result in~\cite[Theorem 4.2]{vDSc}.

In the case when $q_0>0$ and when absorption at $-1$ is certain, which holds for $p\leq1/2$
(see~\cite[Theorem 2.2]{vDSc}),
one gets that if $\eta=1$ then there is no quasi-stationary distribution.
In our model this situation only occurs for $p=1/2$.
Indeed, when $p<1/2$, $p_0+r_0<1$ and $p_0\neq p$, $x=1$ is a~root of the denominator of $m(z)$
in~\eqref{m(z)} if
\begin{gather*}
p=\frac{p_0}{1+p_0-r_0}.
\end{gather*}
But in this case one never has a~mass at $x=1$, since the positive part of the corresponding
residue vanishes in~\eqref{pos part}.
When $p<1/2$, $p_0+r_0<1$ and $p=p_0$, $x=1$ is a~root of the denominator of $m(z)$ when $r_0=p_0$.
In this case the value of the weight in~\eqref{pos part p=p0} is equal to zero.

In case that $\eta<1$, which occurs in our model when $p<1/2$ and $q_0>0$, we have a~one parameter
family of such distributions parametrized by values of $x$ in the interval $\eta<x<1$.

The actual values of the distributions (which depend on~$x$) are given by
\begin{gather*}
\alpha_j(x)=\pi_j(1-x) \frac{Q_j(x)}{q_0},
\end{gather*}
where $Q_j(x)$ are the eigenfunctions given in~\eqref{polynoms}.

\subsection{Ratio limits}

From the very beginning of the work of Karlin and McGregor, involving birth and death processes
with discrete or continuous time parameter, one of the most important advantages of the spectral
approach was the possibility to study
the value of the limit as $n$ grows of quantities such as
\begin{gather}
\label{limitputative}
P^n_{i,j}/P^n_{k,l}.
\end{gather}
As~\cite{KMcG} points out these are dif\/f\/icult quantities to study and back in the 1930's
Doeblin introduced tools to
study the behaviour of the ratio of partial sums of these quantities.
The work in~\cite{KMcG} can thus be seen as getting a~Tauberian theorem to go beyond Doeblin's
results.
In~\cite[Theorem~2]{KMcG} one f\/inds a~proof of the fact that for a~recurrent (and aperiodic) walk
the limit above is given by the ratio $\pi_j/\pi_l$.
In our case, the walk is aperiodic or
non-symmetric when $r_0\neq0$ and recurrent when~\eqref{cond recurrencia} occurs.

Assume that the parameters $r_0,p_0,p, r_0\neq0$, are such that we have at least one nonzero
point mass.
From Sections~\ref{una raiz} and \ref{dos raices} we know that we have just one mass when either
$p=p_0$ and \eqref{una masa p=p0} holds or $p\neq p_0$ along with conditions~\eqref{cond one mass
1} or~\eqref{con one mass 2} and one f\/inds two masses under condition~\eqref{dos masas}.
Then one can use a~result in~\cite[Theorem~3.1]{vDSc1} to see that the limit exists for all values
of the indices $i$, $j$, $k$, $l$ and is equal to
\begin{gather}
\label{limitgen}
\frac{\pi_j}{\pi_l}\frac{Q_i(\eta)Q_j(\eta)}{Q_k(\eta)Q_l(\eta)},
\end{gather}
where $\eta$ is the supremum of all points in the support of the spectral measure, as def\/ined
above.
One can also see that for any values of the parameters $r_0,p_0,p,r_0\neq0$, the limits of the
subsequences resulting by taking $n$ even and odd respectively in~\eqref{limitputative} do exist.
For the symmetric case $r_0=0$, these limits where found in~\cite[Theorem 3]{KMcG} to be equal
to~\eqref{limitgen} for $i-j$ and $k-l$ being simultaneously even or odd.

\subsection{The matrix valued case}

We conclude by calling attention to the last section in~\cite{KMcG} where the authors consider an
extension
of their spectral method to the study of a~walk on the
integers.
They show quite explicitly that the orthogo\-nality measure is to be replaced by a~matrix valued
measure,
which they write down.
This point is noticed, for instance, in the Math.
Reviews note on this paper, MR0100927.
In~\cite{KMcG} one does not f\/ind the corresponding matrix valued orthogo\-nal polynomials or the
extension to the matrix valued case
of the Karlin--McGregor formula.
This was eventually, and independently done much later in~\cite{DRSZ} and \cite{G1}.

\subsection*{Acknowledgements}

The authors thank the anonymous referees for their careful reading and useful remarks.
The work of the f\/irst author is partially supported by MTM2012-36732-C03-03 (Ministerio de
Econom\'ia y Competitividad),
FQM-262, FQM-4643, FQM-7276 (Junta de Andaluc\'ia) and Feder Funds (European
Union) and that of the second author is partially supported
by the Applied Math.
Sciences subprogram of the Of\/f\/ice of Energy
Research, USDOE, under Contract DE--AC03-76SF00098.

\pdfbookmark[1]{References}{ref}
\LastPageEnding


\begin{thebibliography}{99}
\footnotesize\itemsep=0pt

\bibitem{Al1}
Allaway W.R., On f\/inding the distribution function for an orthogonal polynomial
  set, \textit{Pacific~J. Math.} \textbf{49} (1973), 305--310.

\bibitem{Al}
Allaway W.R., The identif\/ication of a~class of orthogonal polynomial sets,
  Ph.D. thesis, University of Alberta, Ann Arbor, MI, 1972.

\bibitem{A1}
Anshelevich M., Bochner--{P}earson-type characterization of the free {M}eixner
  class, \href{http://dx.doi.org/10.1016/j.aam.2010.01.011}{\textit{Adv. in Appl. Math.}} \textbf{46} (2011), 25--45,
  \href{http://arxiv.org/abs/0909.1097}{arXiv:0909.1097}.

\bibitem{A}
Anshelevich M., Free martingale polynomials, \href{http://dx.doi.org/10.1016/S0022-1236(03)00061-2}{\textit{J.~Funct. Anal.}}
  \textbf{201} (2003), 228--261, \href{http://arxiv.org/abs/math.CO/0112194}{math.CO/0112194}.


\bibitem{Apt}
Aptekarev A.I., Asymptotic properties of polynomials orthogonal on a~system of
  contours, and periodic moments of Toda lattices, \href{http://dx.doi.org/10.1070/SM1986v053n01ABEH002918}{\textit{Math. USSR Sb.}}
  \textbf{53} (1986), 233--260.

\bibitem{AW}
Askey R., Wilson J., Some basic hypergeometric orthogonal polynomials that
  generalize {J}acobi polynomials, \textit{Mem. Amer. Math. Soc.} \textbf{54}
  (1985), no.~319, iv+55~pages.

\bibitem{Ch}
Chihara T.S., An introduction to orthogonal polynomials, \textit{Mathematics
  and its Applications}, Vol.~13, Gordon and Breach Science Publishers, New
  York, 1978.

\bibitem{CT}
Cohen J.M., Trenholme A.R., Orthogonal polynomials with a constant recursion
  formula and an application to harmonic analysis, \href{http://dx.doi.org/10.1016/0022-1236(84)90071-5}{\textit{J.~Funct. Anal.}}
  \textbf{59} (1984), 175--184.

\bibitem{DKS}
Damanik D., Killip R., Simon B., Perturbations of orthogonal polynomials with
  periodic recursion coef\/f\/icients, \href{http://dx.doi.org/10.4007/annals.2010.171.1931}{\textit{Ann. of Math.~(2)}} \textbf{171}
  (2010), 1931--2010, \href{http://arxiv.org/abs/math.SP/0702388}{math.SP/0702388}.

\bibitem{DRSZ}
Dette H., Reuther B., Studden W.J., Zygmunt M., Matrix measures and random
  walks with a block tridiagonal transition matrix, \href{http://dx.doi.org/10.1137/050638230}{\textit{SIAM~J. Matrix
  Anal. Appl.}} \textbf{29} (2007), 117--142.

\bibitem{F}
Feller W., An introduction to probability theory and its applications, Vol.~1,
  John Wiley \& Sons, Inc., New York, 1967.

\bibitem{F1}
Feller W., On second order dif\/ferential operators, \textit{Ann. of Math.~(2)}
  \textbf{61} (1955), 90--105.

\bibitem{F2}
Foster F.G., On the stochastic matrices associated with certain queuing
  processes, \href{http://dx.doi.org/10.1214/aoms/1177728976}{\textit{Ann. Math. Stat.}} \textbf{24} (1953), 355--360.

\bibitem{Ge1}
Geronimus J.L., On a set of polynomials, \href{http://dx.doi.org/10.2307/1968164}{\textit{Ann. of Math.~(2)}} \textbf{31}
  (1930), 681--686.

\bibitem{Ge2}
Geronimus J.L., On some equations in f\/inite dif\/ferences and the
  corresponding systems of orthogonal polynomials, \textit{Zap. Mat. Otdel. Fiz.-Mat. Fak. i Kharkov. Mat. Obshch.}
  \textbf{25} (1957), 87--100 (in Russian).

\bibitem{G}
Good I.J., Random motion and analytic continued fractions, \href{http://dx.doi.org/10.1017/S030500410003317X}{\textit{Proc.
  Cambridge Philos. Soc.}} \textbf{54} (1958), 43--47.

\bibitem{G1}
Gr{\"u}nbaum F.A., Random walks and orthogonal polynomials: some challenges, in
  Probability, geometry and integrable systems, \textit{Math. Sci. Res. Inst.
  Publ.}, Vol.~55, Cambridge Univ. Press, Cambridge, 2008, 241--260,
  \href{http://arxiv.org/abs/math.PR/0703375}{math.PR/0703375}.

\bibitem{H}
Harris T.E., First passage and recurrence distributions, \href{http://dx.doi.org/10.1090/S0002-9947-1952-0052057-2}{\textit{Trans. Amer.
  Math. Soc.}} \textbf{73} (1952), 471--486.

\bibitem{HR}
Hodges Jr. J.L., Rosenblatt M., Recurrence-time moments in random walks,
  \textit{Pacific~J. Math.} \textbf{3} (1953), 127--136.

\bibitem{ILMV}
Ismail M.E.H., Masson D.R., Letessier J., Valent G., Birth and death processes
  and orthogonal polynomials, in Orthogonal Polynomials ({C}olumbus, {OH},
  1989), \textit{NATO Adv. Sci. Inst. Ser.~C Math. Phys. Sci.}, Vol.~294,
  Kluwer Acad. Publ., Dordrecht, 1990, 229--255.

\bibitem{K}
Kac M., Random walk and the theory of {B}rownian motion, \textit{Amer. Math.
  Monthly} \textbf{54} (1947), 369--391.

\bibitem{K1}
Karlin S., A f\/irst course in stochastic processes, Academic Press, New York,
  1966.

\bibitem{KMcG}
Karlin S., McGregor J., Random walks, \textit{Illinois~J. Math.} \textbf{3}
  (1959), 66--81.

\bibitem{LR}
Ledermann W., Reuter G.E.H., Spectral theory for the dif\/ferential equations of
  simple birth and death processes, \href{http://dx.doi.org/10.1098/rsta.1954.0001}{\textit{Philos. Trans. Roy. Soc. London.
  Ser.~A.}} \textbf{246} (1954), 321--369.

\bibitem{McK}
McKean Jr. H.P., Elementary solutions for certain parabolic partial
  dif\/ferential equations, \href{http://dx.doi.org/10.1090/S0002-9947-1956-0087012-3}{\textit{Trans. Amer. Math. Soc.}} \textbf{82} (1956),
  519--548.

\bibitem{SY}
Saitoh N., Yoshida H., The inf\/inite divisibility and orthogonal polynomials
  with a constant recursion formula in free probability theory, \textit{Probab.
  Math. Stat.} \textbf{21} (2001), 159--170.

\bibitem{S}
Sodin S., Random matrices, nonbacktracking walks, and orthogonal polynomials,
  \href{http://dx.doi.org/10.1063/1.2819599}{\textit{J.~Math. Phys.}} \textbf{48} (2007), 123503, 21~pages,
  \href{http://arxiv.org/abs/math-ph/0703043}{math-ph/0703043}.

\bibitem{Sz}
Szeg{\"{o}} G., Orthogonal polynomials, \textit{American Mathematical Society,
  Colloquium Publications}, Vol.~23, 4th ed., Amer. Math. Soc., Providence,
  R.I., 1975.

\bibitem{vDSc}
van Doorn E.A., Schrijner P., Geometric ergodicity and quasi-stationarity in
  discrete-time birth-death processes, \href{http://dx.doi.org/10.1017/S0334270000007621}{\textit{J.~Austral. Math. Soc. Ser.~B}}
  \textbf{37} (1995), 121--144.

\bibitem{vDSc1}
van Doorn E.A., Schrijner P., Ratio limits and limiting conditional
  distributions for discrete-time birth-death processes, \href{http://dx.doi.org/10.1006/jmaa.1995.1076}{\textit{J.~Math. Anal.
  Appl.}} \textbf{190} (1995), 263--284.

\end{thebibliography}
\end{document}